\documentclass[12pt]{article}
\usepackage{mathtools} 
\usepackage{amsmath, amssymb} 
\usepackage{amsthm}
\usepackage{enumerate}  
\usepackage{hyperref}
\usepackage{url} 

 \usepackage{tikz}
 \usepackage{tkz-graph}
 \usepackage{caption}

\usepackage{graphicx}

\usepackage{blkarray}

\newcommand\cx{{\mathbb C}}

\newcommand\ints{{\mathbb Z}}
\newcommand\re{{\mathbb R}}

\DeclarePairedDelimiter\abs{\lvert}{\rvert}%
\DeclarePairedDelimiter\norm{\lVert}{\rVert}%

\makeatletter
\let\oldabs\abs
\def\abs{\@ifstar{\oldabs}{\oldabs*}}
\let\oldnorm\norm
\def\norm{\@ifstar{\oldnorm}{\oldnorm*}}
\makeatother

\newcommand\sbs{\subseteq}

\newcommand\comp[1]{{\mkern2mu\overline{\mkern-2mu#1}}}

\newcommand\seq[4]{#1_{#2},#1_{#3},\ldots,#1_{#4}}

\newtheoremstyle{plainsl}%
	{\topsep}
	{\topsep}
	{\slshape} 
	{}
	{\normalfont\bfseries}
	{.}
	{ }
	{}

\swapnumbers

{\theoremstyle{plainsl}
\newtheorem{theorem}{Theorem}[section]
\newtheorem{lemma}[theorem]{Lemma}
\newtheorem{corollary}[theorem]{Corollary}}
{\theoremstyle{remark}
}

\renewcommand\proof{\noindent\textsl{Proof. }}
\newcommand\sqr[2]{{\vbox{\hrule height.#2pt
    \hbox{\vrule width.#2pt height#1pt \kern#1pt
        \vrule width.#2pt}\hrule height.#2pt}}}
\renewcommand\qed{%
	\ifmmode\eqno\sqr53
	\else\nolinebreak\ \hfill\sqr53\medbreak\fi}

\DeclareMathOperator{\rk}{rk}

\DeclareMathOperator{\col}{col}

\newcommand\one{{\bf1}}

\newcommand\grp[1]{\langle #1\rangle}

\title{An Infinite Family of Circulant Graphs with Perfect State Transfer in Discrete Quantum Walks}
\author{Hanmeng Zhan \thanks{Department of Combinatorics and Optimization, University of Waterloo. \texttt{h3zhan@uwaterloo.ca}}}

\begin{document}
\maketitle

\begin{abstract}

We study perfect state transfer in Kendon's model of discrete quantum walks. In particular, we give a characterization of perfect state transfer purely in terms of the graph spectra, and construct an infinite family of $4$-regular circulant graphs that admit perfect state transfer. Prior to our work, the only known infinite families of examples were variants of cycles and diamond chains.
\end{abstract}

\section{Introduction}
Discrete quantum walks were first introduced by Aharonov et al \cite{Aharonov2000} as a quantum analogue of classical random walks. Since then, various models of discrete quantum walks have been proposed \cite{Kendon2003,Szegedy2004,Portugal2015}. In \cite{Bose2002}, Bose posed a scheme for using spin chains to transmit quantum states in quantum computers. Later, Christandl et al \cite{Christandl2004,Christandl2005} studied the problem of perfect state transfer in spin networks with respect to the XY-coupling model. Both continuous and discrete quantum walks were shown to be universal for quantum computation \cite{Childs2009a,Lovett2010,Underwood2010}, for which quantum state transfer plays a role in implementing the universal quantum gate set. In particular, discrete perfect state transfer was first studied by Lovett et al \cite{Lovett2010} on ``wire" structures, where the quantum walk may propagate deterministically.

While there have been numerous results on perfect state transfer in continuous quantum walks \cite{Kay2006,Angeles-Canul2009,Angeles-Canul2010,Kay2011,Cheung2011,Kendon2011,Bachman2012,Coutinho2015,Coutinho2015a,Coutinho2015b}, less is known on the discrete side, as the extra degree of freedom makes it harder to analyze the transition operator. Examples of perfect state transfer in discrete quantum walks have been found on simple structures. In \cite{Kendon2011}, Kendon and Tamon showed that some diamond chains admit perfect state transfer between antipodal vertices. Kurzynski and Wojcik \cite{Kurzynski2011} found perfect state transfer on cycles, and discussed how to convert the position dependence of couplings into the position dependence of coins. Barr et al \cite{Barr2014} investigated discrete quantum walks on variants of cycles, and found some families such as $\comp{K_2}+C_n$ that admit perfect state transfer with appropriately chosen coins and initial states. In \cite{Yalcnkaya2014}, Yalcnkaya and Gedik proposed a scheme to achieve perfect state transfer on paths and cycles using a recovery operator. With various setting of coin flippings, Zhan et al \cite{Zhan2014} also showed that an arbitrary unknown two-qubit state can be perfectly transfered in one-dimensional or two-dimensional lattices. Recently, Stefanak and Skoupy analyzed perfect state transfer on stars \cite{Stefanak2016} and complete bipartite graphs \cite{Stefanak2017} between marked vertices: in $K_{n,n}$, perfect state transfer occurs between any two marked vertices, while in $K_{m,n}$ with $m\ne n$, perfect state transfer only occurs between two marked vertices on the same side.

In this paper, we study perfect state transfer in a simple model proposed by Kendon \cite{Kendon2003}, where each iteration consists of a Grover coin flip followed by an arc-reversal. We give a complete characterization of perfect state transfer in terms of graph spectra. Using this characterization, we construct an infinite family of $4$-regular circulant graphs that admit perfect state transfer. Prior to our work, the only known infinite families of examples were variants of cycles and diamond chains.

\section{Arc-Reversal Grover Walks}
We start by describing Kendon's model of discrete quantum walks \cite{Kendon2003}. While her model applies to all graphs, we will focus on those that are regular. 

Let $X$ be a $d$-regular graph on $n$ vertices. To construct a discrete quantum walk on $X$, we view $X$ as a directed graph, where each edge is replaced by a pair of opposite arcs. A quantum state associated with $X$ is a complex function on its arcs. These states form an inner product vector space, isomorphic to $\cx^n\otimes \cx^d$. Parallel vectors in $\cx^n\otimes \cx^d$ are identified as the same state; we will pick one with unit length as the representative.

A discrete quantum walk is determined by a unitary matrix $U$ acting on $\cx^n\otimes \cx^d$. At step $k$, the system is in state 
\[x_k:=U^kx_0,\]
given initial state $x_0$. We call $U$ the \textsl{transition matrix} of the quantum walk. 

In \cite{Kendon2003}, the transition matrix $U$ is defined to be a product of two sparse unitary matrices. Let $R$ be the permutation matrix that reverses each arc, and $G$ a $d\times d$ unitary matrix of the form
\[G = \frac{2}{d}J - I.\]
The transition matrix of the \textsl{arc-reversal Grover walk} on $X$ is
\[U := R(I\otimes G).\]
We will refer to $R$ as the \textsl{arc-reversal operator},  $G$ as the \textsl{Grover coin}, and $I\otimes G$ the \textsl{coin operator}.

\section{Spectral Decomposition of Two-Reflection Walks}
Notice that both $R$ and $I\otimes G$ have order two. In \cite{Szegedy2004}, Szegedy studied the spectrum of a unitary matrix that is a product of two reflections. Later in Godsil's unpublished notes \cite{Godsil2015a}, he developed some machinery towards finding the spectral decomposition of any matrix lying in the algebra generated by two reflections. We will summarize Godsil's results in this section; a more detailed discussion with proofs will appear in a joint paper by Godsil, Sobchuk and Zhan.

Let $P$ and $Q$ be two projections acting on $\cx^m$. Let $\grp{P, Q}$ denote the matrix algebra generated by $P$ and $Q$. The following well-known result enables us to diagonalize any matrix in $\grp{P,Q}$. 

\begin{lemma}
The vector space $\cx^m$ is a direct sum of $1$- and $2$-dimensional $\grp{P,Q}$-invariant subspaces.
\qed
\end{lemma}

Let $U$ be the product of two reflections about $\col(P)$ and $\col(Q)$, that is, 
\[U:= (2P-I) (2Q-I).\]
To find the spectral decomposition of $U$, we may first decompose $\cx^m$ into a direct sum of $1$- and $2$-dimensional $\grp{P, Q}$-invariant subspaces, and then diagonalize $U$ restricted to them. The $1$-dimensional $\grp{P,Q}$-invariant subspaces are precisely common eigenvectors of $P$ and $Q$. In fact, they span the eigenspaces for $U$ with real eigenvalues, that is, $1$ and $-1$. The $2$-dimensional $\grp{P,Q}$-invariant subspaces provide eigenspaces for $U$ with non-real eigenvalues.

Since $P$ and $Q$ are positive-semidefinite, they have Cholesky decompositions
\[P = KK^*,\quad Q=LL^*,\]
for some matrices $K$ and $L$ with orthonormal columns. Define
\[S:= L^*K.\]
This matrix largely determines the spectrum of $U$.

\begin{lemma}\label{Lem_1es}
The $1$-eigenspace of $U$ is the direct sum
\[(\col(P)\cap\col(Q))\oplus(\ker(P)\cap\ker(Q)),\]
which has dimension 
\[m-\rk(P)-\rk(Q)+2\dim(\col(P)\cap \col(Q)).\]
Moreover, the map $y\mapsto Ly$ is an isomorphism from the $1$-eigenspace of $SS^*$ to $\col(P)\cap\col(Q)$.
\qed
\end{lemma}

\begin{lemma}\label{Lem_-1es}
The $(-1)$-eigenspace of $U$ is the direct sum
\[(\col(P)\cap\ker(Q))\oplus(\ker(P)\cap\col(Q)),\]
which has dimension 
\[\rk(P)+\rk(Q) - 2\rk(S).\]
Moreover, $y\mapsto Ky$ is an isomorphism from $\ker(S)$ to $\col(P)\cap\ker(Q)$, and $y\mapsto L^*y$ is an isomorphism from $\ker(S^*)$ to $\ker(P)\cap\col(Q)$.
\qed
\end{lemma}

\begin{lemma} \label{Lem_nonreal}
The dimensions of the eigenspaces for $U$ with non-real eigenvalues sum to 
\[2\rk(S)-2\dim(\col(P)\cap\col(Q)).\]
Let $\mu\in (0,1)$ be an eigenvalue of $SS^*$. Let $\theta$ be such that $\cos(\theta)=2\mu-1$. The map 
\[y\mapsto ((\cos(\theta)+1)I - (e^{i\theta}+1)P)Ly\]
is an isomorphism from the $\mu$-eigenspace of $SS^*$ to the $e^{i\theta}$-eigenspace of $U$, and the map
\[y\mapsto ((\cos(\theta)+1)I - (e^{-i\theta}+1)P)Ly\]
is an isomorphism from the $\mu$-eigenspace of $SS^*$ to the $e^{-i\theta}$-eigenspace of $U$. 
\qed
\end{lemma}

\section{Graph Spectra vs Walk Spectra}
We apply the previous results to the transition matrix of an arc-reversal Grover walk:
\[U = R(I\otimes G).\]
While $U$ is not an obvious function in the adjacency matrix of the graph $X$, we find a spectral correspondence between $U$ and $X$. More specifically, eigenvalues of $X$ provide the real parts of eigenvalues of $U$, and eigenvectors of $X$ can be lifted to eigenvectors of $U$ by two incidence matrices. 

To start, we introduce four incidence matrices: the tail-arc incidence matrix $D_t$, the head-arc incidence matrix $D_h$, the arc-edge incidence matrix $M$, and the vertex-edge incidence matrix $B$. 

The tail-arc incidence matrix $D_t$, and the head-arc incidence matrix $D_h$, are two matrices with rows indexed by the vertices, and columns by the arcs. If $u$ is a vertex and $e$ is an edge, then $(D_t)_{u,e}=1$ if $u$ is the initial vertex of $e$, and $(D_h)_{u,e}=1$ if $e$ ends on $u$.

The arc-edge incidence matrix $M$ is a matrix with rows indexed by the arcs and columns by the edges. If $a$ is an arc and $e$ is an edge, then $M_{a,e}=1$ if $a$ is one direction of $e$.

The vertex-edge incidence matrix $B$ is a matrix with rows indexed by the vertices and columns by the edges. If $u$ is a vertex and $e$ is an edge, then $B_{u,e}=1$ if $u$ is one endpoints of $e$.

As an example, the following are the four incidence matrices associated with $K_3$ with vertices $\{0,1,2\}$. 

\[D_t=
\begin{blockarray}{ccccccc}
& (0,1) & (0,2) & (1,0) & (1,2) & (2,0) & (2,1) \\
\begin{block}{c(cccccc)}
0 & 1 & 1 & 0 & 0 & 0 & 0\\
1 & 0 & 0 & 1 & 1 & 0 & 0\\
2 & 0 & 0 & 0 & 0 & 1 & 1\\
\end{block}
\end{blockarray}\]
\[D_h=
\begin{blockarray}{ccccccc}
& (0,1) & (0,2) & (1,0) & (1,2) & (2,0) & (2,1) \\
\begin{block}{c(cccccc)}
0 & 0 & 0 & 1 & 0 & 1 & 0\\
1 & 1 & 0 & 0 & 0 & 0 & 1\\
2 & 0 & 1 & 0 & 1 & 0 & 0\\
\end{block}
\end{blockarray}\]
\[M = 
\begin{blockarray}{cccc}
& \{0,1\} & \{0,2\} & \{1,2\}\\
\begin{block}{c(ccc)}
(0,1) & 1 & 0 & 0\\
(0,2) & 0 & 1 & 0\\
(1,0) & 1 & 0 & 0\\
(1,2) & 0 & 0 & 1\\
(2,0) & 0 & 1 & 0\\
(2,1) & 0 & 0 & 1\\
\end{block}
\end{blockarray}\]
\[B =
\begin{blockarray}{cccc}
& \{0,1\} & \{0,2\} & \{1,2\}\\
\begin{block}{c(ccc)}
0 & 1 & 1 & 0\\
1 & 1 & 0 & 1\\
2 & 0 & 1 & 1\\
\end{block}
\end{blockarray}
\]

Next, we list some useful identities about these incidence matrices. Let $A$ be the adjacency matrix of $X$. One can verify the following. 

\begin{lemma}\label{Lem_ids_incs}
	We have
	\begin{enumerate}[(i)]
		\item $D_t^TD_t = D_h^TD_h=dI$
		\item $M^TM = 2I$,
		\item $D_t D_h^T = D_hD_t^T = A$,
		\item $BB^T = A + dI$
		\item $D_tM = D_hM = B$,
		\item $D_tR = D_h$, 
		\item $R=MM^T - I$, and
		\item $I \otimes G = \frac{2}{d} D_t^TD_t- I$.
		\qed
	\end{enumerate}
\end{lemma}

As a consequence, $R$ is a reflection about $\col(M)$, while $I\otimes G$ is a reflection about $\col(D_t^T)$. We now prove the spectral relation between $U$ and $X$. The following theorem shows that all eigenspaces of $U$ with non-real eigenvalues are completely determined by the eigenspaces of $X$ with eigenvalues in $(-d,d)$. It also gives a concrete description on how to ``lift" eigenvalues and eigenvectors of $X$ to those of $U$.

\begin{theorem}
Let $y$ be an eigenvector for $X$ with eigenvalue $\lambda\in(-d,d)$. Suppose $\lambda = d\cos(\theta)$ for some $\theta\in\re$. Then 
\[D_t^Ty - e^{i\theta}D_h^Ty\]
is an eigenvector for $U$ with eigenvalue $e^{i\theta}$, and 
\[D_t^Ty - e^{-i\theta}D_h^Ty\]
is an eigenvector for $U$ with eigenvalue $e^{-i\theta}$.
\end{theorem}
\proof
Let 
\[K:=\frac{1}{\sqrt{2}}M,\quad L:=\frac{1}{\sqrt{d}}D_t^T,\quad S:=L^*K.\]
According to Lemma \ref{Lem_nonreal}, eigenspace for $U$ with non-real eigenvalues are determined by eigenspaces for
\[SS^*= \frac{1}{2d}BB^T = \frac{1}{2d}(A+dI).\]
Let
\[\mu:=\frac{\lambda+d}{2d}.\]
Then $0<\mu<1$ and $2\mu-1 = \cos(\theta)$. Moreover, 
\[Ay = \lambda y\]
if and only if 
\[SS^*y = \mu y.\]
Thus, using identities in Lemma \ref{Lem_ids_incs}, we obtain two eigenvectors for $U$ as stated.
\qed

After normalization, we obtain the eigenprojections for non-real eigenvalues of $U$.

\begin{corollary}\label{Lem_arcrev_nonreal}
Let $\lambda$ be an eigenvalue of $X$ that is neither $d$ nor $-d$. Let $E_{\lambda}$ be the orthogonal projection onto the $\lambda$-eigenspace of $X$. Suppose $\lambda=d\cos(\theta)$ for some $\theta\in\re$. Then the $e^{i\theta}$-eigenprojection of $U$ is
\[\frac{1}{2d\sin^2(\theta)}(D_t-e^{i\theta}D_h)^TE_{\lambda}(D_t-e^{-i\theta}D_h),\]
and the $e^{-i\theta}$-eigenprojection of $U$ is
\[\frac{1}{2d\sin^2(\theta)}(D_t-e^{-i\theta}D_h)^TE_{\lambda}(D_t-e^{i\theta}D_h).\tag*{\sqr53}\]
\end{corollary}

We also characterize the $(\pm 1)$-eigenspaces of $U$. In particular, their multiplicities depends on parameters of $X$. 

\begin{lemma}\label{Lem_arcrev_1es}
The $1$-eigenspace of $U$ is 
\[(\col(M)\cap \col(D_t^T)) \oplus (\ker(M^T) \cap \ker(D_t))\]
with dimension
\[\frac{nd}{2} - n + 2.\]
Moreover, the projection onto $\col(M)\cap \col(D_t^T)$ is given by
\[\frac{1}{d}D_t^TE_d D_t = \frac{1}{nd}J,\]
where $E_d$ is the projection onto the $d$-eigenspace of $X$. 
\end{lemma}
\proof
By Lemma \ref{Lem_1es}, the $1$-eigenspace is the direct sum:
\[(\col(M)\cap \col(D_t^T)) \oplus (\ker(M^T) \cap \ker(D_t)),\]
where
\[\col(M)\cap \col(D_t^T) = D_t \col(E_d).\]
Note that $\col(D_t^T)$ consists of vectors that are constant over the outgoing arcs of each vertex, and $\col(M)$ consists of vectors that are constant over each pair of opposite arcs. Since $X$ is connected,
\[\col(M)\cap \col(D_t^T) = \mathrm{span}\{\one\}.\]
The multiplicity follows from the fact that $\rk(M) = nd/2$ and $\rk(D_t) = n$.
\qed

\begin{lemma}\label{Lem_arcrev_-1es}
If $X$ is bipartite, the $(-1)$-eigenspace of $U$ is 
\[M\ker(B)\oplus D_t^T\ker(B^T)\]
with dimension 
\[\frac{nd}{2} - n +2.\]
Moreover, the projection onto $ D_t^T\ker(B^T)$ is given by
\[\frac{1}{d}D_t^T E_{-d}D_t,\]
where $E_{-d}$ is the projection onto the $(-d)$-eigenspace of $X$. If $X$ is not bipartite, the $(-1)$-eigenspace of $U$ is
\[M\ker(B),\]
with dimension
\[\frac{nd}{2}-n.\]
\end{lemma}
\proof
By Lemma \ref{Lem_-1es}, the $(-1)$-eigenspace of $U$ is 
\[M\ker(B)\oplus D_t^T\ker(B^T),\]
where
\[\ker(B^T) = \col(E_{-d}).\]
Note that $\rk(B)=n-1$ if $X$ is bipartite, and $\rk(B)=n$ otherwise. 
\qed

\section{Perfect State Transfer}
In this section, we derive necessary and sufficient conditions for perfect state transfer to occur on arc-reversal Grover walks. The techniques we use are very similar to those employed in continuous quantum walks. For a thorough treatment of continuous-time perfect state transfer, see Coutinho's Ph.D. thesis \cite{Coutinho2014}. 

Suppose we start with a state that ``concentrates on" $u$, that is, a complex function that sends all but the outgoing arcs of $u$ to zero. In theory, this state could be $e_u\otimes x$ for any unit vector $x$. However, it is more practical to prepare a uniform superposition over the outgoing arcs of $u$:
\[\frac{1}{\sqrt{d}} e_u \otimes \one.\]
Formally, if there is a unit vector $x\in \cx^d$ such that
\[U^k \left(\frac{1}{\sqrt{d}}e_u\otimes \one\right) = e_v\otimes x,\]
then we say $X$ admits \textsl{perfect state transfer} from $u$ to $v$ if $u\ne v$, and $X$ is \textsl{periodic} at $u$ if $u=v$. While this definition does not impose further condition on the final state,  in the arc-reversal Grover walk, the only possible choice of $x$ is 
\[\frac{1}{\sqrt{d}}\one,\]
as we show now.

\begin{lemma}
If $X$ admits perfect state transfer from $u$ to $v$ at time $k$, then
\[U^k \left(\frac{1}{\sqrt{d}}e_u\otimes \one\right) = \left(\frac{1}{\sqrt{d}}e_v \otimes \one\right).\]
\end{lemma}
\proof
Suppose 
\[U^k \left(\frac{1}{\sqrt{d}}e_u\otimes \one\right) = e_v\otimes x.\]
Since $U$ has real entries, all entries in $x$ are also real. Moreover, as $\one\otimes \one$ is an eigenvector for $U$ with eigenvalue $1$, 
\[\left<\one\otimes \one ,\frac{1}{\sqrt{d}} e_u\otimes \one\right> =\left<\one \otimes \one, U^k\left(\frac{1}{\sqrt{d}} e_u\otimes \one\right)\right> =  \grp{\one\otimes \one, e_v\otimes x}.\]
It follows that
\[\grp{\one, x} =\sqrt{d}.\]
On the other hand, by Cauchy-Schwarz, 
\[ \abs{\grp{\one ,x}}\le \norm{\one} \norm{x} = \sqrt{d},\]
with equality held if and only if $x$ is a scalar multiple of $\one$. Therefore $x$ must be equal to $\one$.
\qed

Notice that both the initial state and the final state lie in $\col(D_t^T)$, so an equivalent definition for perfect state transfer from $u$ to $v$ at time $k$ is
\[U^k D_t^T e_u = D_t^T e_v.\]
Our characterization of perfect state transfer relies heavily on this observation.

\begin{lemma}\label{Lem_escollapse}
Let $\lambda = d\cos(\theta)$ be an eigenvalue of $X$ that is neither $d$ nor $-d$. Let $E_{\lambda}$ be the projection onto the $\lambda$-eigenspace of $X$, and let $F_{\pm}$ be the projection onto the $e^{\pm i\theta}$-eigenspace of $U$. Then 
\[D_tF_{\pm} D_t^T = \frac{d}{2} E_{\lambda}.\]
\end{lemma}
\proof
By Lemma \ref{Lem_arcrev_nonreal},
\begin{align*}
2d\sin^2(\theta) D_t F_{\pm} D_t^T&=D_t (D_t-e^{\pm i\theta})^TE_{\lambda} (D_t-e^{\mp i\theta}D_h)\\
&=(dI - e^{\pm i\theta} A) E_{\lambda} (dI - e^{\mp i\theta} A)\\
&=d^2\abs{1-e^{i\theta}\cos(\theta)}^2 E_{\lambda}\\
&=d^2\sin^2(\theta).
\tag*{\sqr53}
\end{align*}

\begin{theorem}\label{Thm_pst}
Suppose the spectral decomposition of $X$ is
\[A = \sum_{\lambda} \lambda E_{\lambda}.\]
Then $X$ admit perfect state transfer from $u$ to $v$ at time $k$ if and only if the following hold.
\begin{enumerate}[(i)]
\item For each $\lambda$, we have $E_{\lambda} e_u = \pm E_{\lambda} e_v$.
\item If $E_{\lambda} e_u = E_{\lambda} e_v\ne 0$, then $\lambda =d\cos(j\pi/k)$ for some even integer $j$.
\item If $E_{\lambda} e_u = -E_{\lambda} e_v\ne 0$, then $\lambda = d\cos(j\pi /k)$ for some odd integer $j$.
\end{enumerate}
\end{theorem}
\proof
Consider the spectral decomposition of $U$:
\[U = \sum_r e^{i\theta_r} F_r.\]
There is perfect state transfer from $u$ to $v$ at time $k$ if and only if 
\[\sum_r e^{ik\theta_r} F_r D_t^Te_u = D_t^Te_v,\]
or equivalently, for each $r$,
\begin{equation}
e^{ik\theta_r} F_r D_t^Te_u = F_r D_t^Te_v.\label{Eqn_sp}
\end{equation}

Suppose $e^{i\theta_r}=1$. Equation (\ref{Eqn_sp}) says that
\[F_r D_t^Te_u = F_r D_t^Te_v.\]
By Lemma \ref{Lem_arcrev_1es}, this holds if and only if 
\[\frac{1}{nd}Je_u = D_t^TE_de_u = D_t^TE_d e_v= \frac{1}{nd}Je_v \ne 0,\]
if and only if 
\[E_d e_u = E_d e_v\ne 0.\]
Clearly $d=\cos(0)$, which satisfies (ii).

Suppose $e^{i\theta_r}=-1$. By Lemma \ref{Lem_arcrev_-1es}, 
\[F_r D_t^T= \frac{1}{d}D_t^TE_{-d}D_t D_t^T = D_t^TE_d.\]
Thus Equation (\ref{Eqn_sp}) holds if and only if 
\[(-1)^k F_r D_t^T e_u = F_r D_t^T e_v,\]
that is,
\[(-1)^k D_t^TE_{-d}e_u = D_t^T E_{-d}e_v.\]
If $X$ is not bipartite, then $E_{-d}=0$ and 
\[F_r D_t^Te_u = F_r D_t^T e_v=0.\]
Otherwise, 
\[E_{-d}e_u = E_{-d}e_v \ne 0\]
if $u$ and $v$ are in the same color class, and
\[E_{-d}e_u = -E_{-d}e_v \ne 0\]
if they are in different color classes. Clearly 
\[-d = \cos\left(\frac{k\pi}{k}\right),\]
which satisfies (i) and (ii).

Finally suppose $e^{i\theta_r}\ne \pm 1$. Equation (\ref{Eqn_sp}) says that
\[e^{ik\theta_r} F_r D_t^Te_u = F_r D_t^Te_v.\]
By Lemma \ref{Lem_escollapse},
\[D_tF_r D_t^T = \frac{d}{2} E_{\lambda},\]
so
\[\frac{de^{ik\theta_r}}{2} (E_{\lambda})_{uu}=e^{ik\theta_r}\left<F_rD_t^Te_u, D_t^Te_u\right> = \left<F_rD_t^Te_v, D_t^Te_u\right>=\frac{d}{2}(E_{\lambda})_{uv} \in \re.\]
Therefore Equation (\ref{Eqn_sp}) holds if and only if one of the following occurs:
\begin{enumerate}[(a)]
\item $E_{\lambda}e_u = E_{\lambda}e_v=0$;
\item $E_{\lambda}e_u = E_{\lambda}e_v\ne 0$, and $e^{ik\theta_r}=1$;
\item $E_{\lambda}e_u = -E_{\lambda} e_v\ne 0$, and $e^{ik\theta_r}=-1$.
\qed
\end{enumerate}

The three conditions in Theorem \ref{Thm_pst} are symmetric in $u$ and $v$. As a consequence, perfect state transfer is symmetric in the initial and final state, and it implies periodicity at both vertices.

\begin{corollary}
If there is perfect state transfer from $u$ to $v$ at time $k$, then there is perfect state transfer from $v$ to $u$ at time $k$, and $X$ is periodic at both $u$ and $v$ at time $2k$.
\qed
\end{corollary} 

Let $X$ be a graph with spectral decomposition
\[A = \sum_{\lambda} \lambda E_{\lambda}.\]
The \textsl{eigenvalue support} of a vertex $u$, defined by Godsil \cite{Godsil2011a}, is the set
\[\{\lambda: E_{\lambda} e_u\ne 0\}.\]
Let $\phi(t)$ be the characteristic polynomial of $X$, and $\phi_u(t)$ the characteristic polynomial of the vertex-deleted subgraph $X-u$. It is shown in \cite{Godsil2001} that the eigenvalue support of $u$ consists of roots of the following polynomial:
\[\psi_u(t):=\frac{\phi(t)}{\gcd(\phi(t), \phi_u(t))}.\]
Thus, Theorem \ref{Thm_pst} gives necessary and sufficient conditions on $\psi_u(t)$ for $X$ to be periodic at $u$.

\begin{theorem}\label{Thm_evsp}
Suppose $\psi_u(t)$ has degree $\ell$. Then vertex $u$ is periodic at time $k$ if and only if the polynomial
\[z^{\ell}\psi_u\left(\frac{d}{2} \left(z + \frac{1}{z}\right)\right)\]
is a factor of $z^k-1$.
\end{theorem}
\proof
Setting $u=v$ in Theorem \ref{Thm_pst}, we see that $u$ is periodic at time $k$ if and only if each eigenvalue $\lambda$ in the eigenvalue support of $u$ is of the form
\[\lambda = \frac{d}{2} (e^{j\pi i/k} + e^{-j\pi i/k}),\]
for some even integer $j$, or equivalently, 
\[z^{\ell}\psi_u\left(\frac{d}{2} \left(z + \frac{1}{z}\right)\right)\]
divides $z^k-1$.
\qed

We say two vertices $u$ and $v$ are \textsl{strongly cospectral} if 
\[E_{\lambda} e_u = \pm E_{\lambda} e_v\]
for each eigenvalue $\lambda$ of $X$. Strongly cospectrality has been thoroughly studied by Godsil and Smith \cite{Godsil2017Smith}; we cite a useful characterization below.

\begin{theorem}
Two vertices $u$ and $v$ are strongly cospectral if and only if the following hold.
\begin{enumerate}[(i)]
\item $u$ and $v$ are cospectral.
\item For every eigenvalue $\lambda$ of $X$, the vectors $E_{\lambda}e_u$ and $E_{\lambda}e_v$ are parallel.
\qed
\end{enumerate}
\end{theorem}

\section{An Infinite Family}

Our characterization of perfect state transfer leads us to consider regular graphs whose eigenvalues are given by real parts of $2k$-th roots of unity. A \textsl{circulant graph} $X=X(\ints_n, \{\seq{c}{1}{2}{d}\})$ is a Cayley graph over $\ints_n$ with connection set 
\[\{\seq{c}{1}{2}{d}\}\sbs \ints_n.\]
If $\psi$ is a character of $\ints_n$, then $\psi$ is also an eigenvector for $X$ with eigenvalue
\[\psi(c_1)+\cdots+\psi(c_d).\] 
Note that this is a sum of real parts of $n$-th roots of unity. We show that circulant graphs whose connection sets satisfy a simple condition admit perfect state transfer. 

\begin{theorem}\label{Thm_pstcirc}
	Let $\ell$ be an odd integer. For any distinct integers $a$ and $b$ such that $a+b=\ell$, the circulant graph $X(\ints_{2\ell}, \pm\{a,b\})$ admits perfect state transfer at time $2\ell$ from vertex $0$ to vertex $\ell$.
\end{theorem}
\proof
The eigenvalues of $X$ are
\begin{align*}
\lambda_j &= e^{aj\pi/\ell} +e^{-aj\pi/\ell} + e^{bj\pi/\ell} +e^{-bj\pi/\ell} \\
&= 2\cos\left(\frac{aj\pi}{\ell}\right) + 2\cos\left(\frac{bj\pi}{\ell}\right),
\end{align*}
for $j=0,1,\cdots,2n-1$. Since $\ell$ is odd and $a+b=\ell$, when $j$ is odd, 
\[\lambda_j=0=4\cos\left(\frac{\ell \pi}{2\ell}\right),\] and when $j$ is even, 
\[\lambda_j = 4\cos\left(\frac{2aj\pi}{2\ell}\right).\]
It suffices to check the parity condition in Theorem \ref{Thm_pst} for each eigenvector of $X$. Since $a+b=\ell$, vertex $u$ and $u+\ell$ have the same neighbors, so
\[A(e_u-e_{u+\ell})=0.\]
We see from the multiplicity of $0$ that for $u=0,1,\cdots, \ell-1$, the vectors $e_u - e_{u+\ell}$ form an orthogonal basis for $\ker(A)$. Thus $y_u=-y_v$ if $y$ is an eigenvector for $X$ with eigenvalue $0$, and $y_u=y_v$ if $y$ is any other eigenvector  for $X$.
\qed

Figures \ref{X6_12}, \ref{X10_14}, \ref{X10_23}, \ref{X14_16}, \ref{X14_25}, and \ref{X14_34} are the first few examples in the infinite family with perfect state transfer between the marked vertices.

\begin{figure}
\centering
\begin{minipage}{0.5\textwidth}
\centering
\begin{tikzpicture}
\definecolor{cv0}{rgb}{0.0,0.0,0.0}
\definecolor{cfv0}{rgb}{1.0,1.0,1.0}
\definecolor{clv0}{rgb}{0.0,0.0,0.0}
\definecolor{cv1}{rgb}{0.0,0.0,0.0}
\definecolor{cfv1}{rgb}{1.0,1.0,1.0}
\definecolor{clv1}{rgb}{0.0,0.0,0.0}
\definecolor{cv2}{rgb}{0.0,0.0,0.0}
\definecolor{cfv2}{rgb}{1.0,1.0,1.0}
\definecolor{clv2}{rgb}{0.0,0.0,0.0}
\definecolor{cv3}{rgb}{0.0,0.0,0.0}
\definecolor{cfv3}{rgb}{1.0,1.0,1.0}
\definecolor{clv3}{rgb}{0.0,0.0,0.0}
\definecolor{cv4}{rgb}{0.0,0.0,0.0}
\definecolor{cfv4}{rgb}{1.0,1.0,1.0}
\definecolor{clv4}{rgb}{0.0,0.0,0.0}
\definecolor{cv5}{rgb}{0.0,0.0,0.0}
\definecolor{cfv5}{rgb}{1.0,1.0,1.0}
\definecolor{clv5}{rgb}{0.0,0.0,0.0}
\definecolor{cv0v1}{rgb}{0.0,0.0,0.0}
\definecolor{cv0v2}{rgb}{0.0,0.0,0.0}
\definecolor{cv0v4}{rgb}{0.0,0.0,0.0}
\definecolor{cv0v5}{rgb}{0.0,0.0,0.0}
\definecolor{cv1v2}{rgb}{0.0,0.0,0.0}
\definecolor{cv1v3}{rgb}{0.0,0.0,0.0}
\definecolor{cv1v5}{rgb}{0.0,0.0,0.0}
\definecolor{cv2v3}{rgb}{0.0,0.0,0.0}
\definecolor{cv2v4}{rgb}{0.0,0.0,0.0}
\definecolor{cv3v4}{rgb}{0.0,0.0,0.0}
\definecolor{cv3v5}{rgb}{0.0,0.0,0.0}
\definecolor{cv4v5}{rgb}{0.0,0.0,0.0}
\Vertex[style={minimum size=1.0cm,draw=cv0,fill=cfv0,text=clv0,shape=circle},LabelOut=false,L=\hbox{$0$},x=5.0cm,y=2.5cm]{v0}
\Vertex[style={minimum size=1.0cm,draw=cv1,fill=cfv1,text=clv1,shape=circle},LabelOut=false,L=\hbox{$1$},x=3.75cm,y=5.0cm]{v1}
\Vertex[style={minimum size=1.0cm,draw=cv2,fill=cfv2,text=clv2,shape=circle},LabelOut=false,L=\hbox{$2$},x=1.25cm,y=5.0cm]{v2}
\Vertex[style={minimum size=1.0cm,draw=cv3,fill=cfv3,text=clv3,shape=circle},LabelOut=false,L=\hbox{$3$},x=0.0cm,y=2.5cm]{v3}
\Vertex[style={minimum size=1.0cm,draw=cv4,fill=cfv4,text=clv4,shape=circle},LabelOut=false,L=\hbox{$4$},x=1.25cm,y=0.0cm]{v4}
\Vertex[style={minimum size=1.0cm,draw=cv5,fill=cfv5,text=clv5,shape=circle},LabelOut=false,L=\hbox{$5$},x=3.75cm,y=0.0cm]{v5}
\Edge[lw=0.02cm,style={color=cv0v1,},](v0)(v1)
\Edge[lw=0.02cm,style={color=cv0v2,},](v0)(v2)
\Edge[lw=0.02cm,style={color=cv0v4,},](v0)(v4)
\Edge[lw=0.02cm,style={color=cv0v5,},](v0)(v5)
\Edge[lw=0.02cm,style={color=cv1v2,},](v1)(v2)
\Edge[lw=0.02cm,style={color=cv1v3,},](v1)(v3)
\Edge[lw=0.02cm,style={color=cv1v5,},](v1)(v5)
\Edge[lw=0.02cm,style={color=cv2v3,},](v2)(v3)
\Edge[lw=0.02cm,style={color=cv2v4,},](v2)(v4)
\Edge[lw=0.02cm,style={color=cv3v4,},](v3)(v4)
\Edge[lw=0.02cm,style={color=cv3v5,},](v3)(v5)
\Edge[lw=0.02cm,style={color=cv4v5,},](v4)(v5)
\end{tikzpicture}
	\captionof{figure}{$X(\ints_6, \{1,2,4,5\})$}
	\label{X6_12}
\end{minipage}%
\begin{minipage}{0.5\textwidth}
\centering
\begin{tikzpicture}
\definecolor{cv0}{rgb}{0.0,0.0,0.0}
\definecolor{cfv0}{rgb}{1.0,1.0,1.0}
\definecolor{clv0}{rgb}{0.0,0.0,0.0}
\definecolor{cv1}{rgb}{0.0,0.0,0.0}
\definecolor{cfv1}{rgb}{1.0,1.0,1.0}
\definecolor{clv1}{rgb}{0.0,0.0,0.0}
\definecolor{cv2}{rgb}{0.0,0.0,0.0}
\definecolor{cfv2}{rgb}{1.0,1.0,1.0}
\definecolor{clv2}{rgb}{0.0,0.0,0.0}
\definecolor{cv3}{rgb}{0.0,0.0,0.0}
\definecolor{cfv3}{rgb}{1.0,1.0,1.0}
\definecolor{clv3}{rgb}{0.0,0.0,0.0}
\definecolor{cv4}{rgb}{0.0,0.0,0.0}
\definecolor{cfv4}{rgb}{1.0,1.0,1.0}
\definecolor{clv4}{rgb}{0.0,0.0,0.0}
\definecolor{cv5}{rgb}{0.0,0.0,0.0}
\definecolor{cfv5}{rgb}{1.0,1.0,1.0}
\definecolor{clv5}{rgb}{0.0,0.0,0.0}
\definecolor{cv6}{rgb}{0.0,0.0,0.0}
\definecolor{cfv6}{rgb}{1.0,1.0,1.0}
\definecolor{clv6}{rgb}{0.0,0.0,0.0}
\definecolor{cv7}{rgb}{0.0,0.0,0.0}
\definecolor{cfv7}{rgb}{1.0,1.0,1.0}
\definecolor{clv7}{rgb}{0.0,0.0,0.0}
\definecolor{cv8}{rgb}{0.0,0.0,0.0}
\definecolor{cfv8}{rgb}{1.0,1.0,1.0}
\definecolor{clv8}{rgb}{0.0,0.0,0.0}
\definecolor{cv9}{rgb}{0.0,0.0,0.0}
\definecolor{cfv9}{rgb}{1.0,1.0,1.0}
\definecolor{clv9}{rgb}{0.0,0.0,0.0}
\definecolor{cv0v1}{rgb}{0.0,0.0,0.0}
\definecolor{cv0v4}{rgb}{0.0,0.0,0.0}
\definecolor{cv0v6}{rgb}{0.0,0.0,0.0}
\definecolor{cv0v9}{rgb}{0.0,0.0,0.0}
\definecolor{cv1v2}{rgb}{0.0,0.0,0.0}
\definecolor{cv1v5}{rgb}{0.0,0.0,0.0}
\definecolor{cv1v7}{rgb}{0.0,0.0,0.0}
\definecolor{cv2v3}{rgb}{0.0,0.0,0.0}
\definecolor{cv2v6}{rgb}{0.0,0.0,0.0}
\definecolor{cv2v8}{rgb}{0.0,0.0,0.0}
\definecolor{cv3v4}{rgb}{0.0,0.0,0.0}
\definecolor{cv3v7}{rgb}{0.0,0.0,0.0}
\definecolor{cv3v9}{rgb}{0.0,0.0,0.0}
\definecolor{cv4v5}{rgb}{0.0,0.0,0.0}
\definecolor{cv4v8}{rgb}{0.0,0.0,0.0}
\definecolor{cv5v6}{rgb}{0.0,0.0,0.0}
\definecolor{cv5v9}{rgb}{0.0,0.0,0.0}
\definecolor{cv6v7}{rgb}{0.0,0.0,0.0}
\definecolor{cv7v8}{rgb}{0.0,0.0,0.0}
\definecolor{cv8v9}{rgb}{0.0,0.0,0.0}
\Vertex[style={minimum size=1.0cm,draw=cv0,fill=cfv0,text=clv0,shape=circle},LabelOut=false,L=\hbox{$0$},x=5.0cm,y=2.5cm]{v0}
\Vertex[style={minimum size=1.0cm,draw=cv1,fill=cfv1,text=clv1,shape=circle},LabelOut=false,L=\hbox{$1$},x=4.5225cm,y=4.0451cm]{v1}
\Vertex[style={minimum size=1.0cm,draw=cv2,fill=cfv2,text=clv2,shape=circle},LabelOut=false,L=\hbox{$2$},x=3.2725cm,y=5.0cm]{v2}
\Vertex[style={minimum size=1.0cm,draw=cv3,fill=cfv3,text=clv3,shape=circle},LabelOut=false,L=\hbox{$3$},x=1.7275cm,y=5.0cm]{v3}
\Vertex[style={minimum size=1.0cm,draw=cv4,fill=cfv4,text=clv4,shape=circle},LabelOut=false,L=\hbox{$4$},x=0.4775cm,y=4.0451cm]{v4}
\Vertex[style={minimum size=1.0cm,draw=cv5,fill=cfv5,text=clv5,shape=circle},LabelOut=false,L=\hbox{$5$},x=0.0cm,y=2.5cm]{v5}
\Vertex[style={minimum size=1.0cm,draw=cv6,fill=cfv6,text=clv6,shape=circle},LabelOut=false,L=\hbox{$6$},x=0.4775cm,y=0.9549cm]{v6}
\Vertex[style={minimum size=1.0cm,draw=cv7,fill=cfv7,text=clv7,shape=circle},LabelOut=false,L=\hbox{$7$},x=1.7275cm,y=0.0cm]{v7}
\Vertex[style={minimum size=1.0cm,draw=cv8,fill=cfv8,text=clv8,shape=circle},LabelOut=false,L=\hbox{$8$},x=3.2725cm,y=0.0cm]{v8}
\Vertex[style={minimum size=1.0cm,draw=cv9,fill=cfv9,text=clv9,shape=circle},LabelOut=false,L=\hbox{$9$},x=4.5225cm,y=0.9549cm]{v9}
\Edge[lw=0.02cm,style={color=cv0v1,},](v0)(v1)
\Edge[lw=0.02cm,style={color=cv0v4,},](v0)(v4)
\Edge[lw=0.02cm,style={color=cv0v6,},](v0)(v6)
\Edge[lw=0.02cm,style={color=cv0v9,},](v0)(v9)
\Edge[lw=0.02cm,style={color=cv1v2,},](v1)(v2)
\Edge[lw=0.02cm,style={color=cv1v5,},](v1)(v5)
\Edge[lw=0.02cm,style={color=cv1v7,},](v1)(v7)
\Edge[lw=0.02cm,style={color=cv2v3,},](v2)(v3)
\Edge[lw=0.02cm,style={color=cv2v6,},](v2)(v6)
\Edge[lw=0.02cm,style={color=cv2v8,},](v2)(v8)
\Edge[lw=0.02cm,style={color=cv3v4,},](v3)(v4)
\Edge[lw=0.02cm,style={color=cv3v7,},](v3)(v7)
\Edge[lw=0.02cm,style={color=cv3v9,},](v3)(v9)
\Edge[lw=0.02cm,style={color=cv4v5,},](v4)(v5)
\Edge[lw=0.02cm,style={color=cv4v8,},](v4)(v8)
\Edge[lw=0.02cm,style={color=cv5v6,},](v5)(v6)
\Edge[lw=0.02cm,style={color=cv5v9,},](v5)(v9)
\Edge[lw=0.02cm,style={color=cv6v7,},](v6)(v7)
\Edge[lw=0.02cm,style={color=cv7v8,},](v7)(v8)
\Edge[lw=0.02cm,style={color=cv8v9,},](v8)(v9)
\end{tikzpicture}
	\captionof{figure}{$X(\ints_{10}, \{1,4,6,9\})$}
	\label{X10_14}
\end{minipage}
\begin{minipage}{0.5\textwidth}
\centering
\begin{tikzpicture}
\definecolor{cv0}{rgb}{0.0,0.0,0.0}
\definecolor{cfv0}{rgb}{1.0,1.0,1.0}
\definecolor{clv0}{rgb}{0.0,0.0,0.0}
\definecolor{cv1}{rgb}{0.0,0.0,0.0}
\definecolor{cfv1}{rgb}{1.0,1.0,1.0}
\definecolor{clv1}{rgb}{0.0,0.0,0.0}
\definecolor{cv2}{rgb}{0.0,0.0,0.0}
\definecolor{cfv2}{rgb}{1.0,1.0,1.0}
\definecolor{clv2}{rgb}{0.0,0.0,0.0}
\definecolor{cv3}{rgb}{0.0,0.0,0.0}
\definecolor{cfv3}{rgb}{1.0,1.0,1.0}
\definecolor{clv3}{rgb}{0.0,0.0,0.0}
\definecolor{cv4}{rgb}{0.0,0.0,0.0}
\definecolor{cfv4}{rgb}{1.0,1.0,1.0}
\definecolor{clv4}{rgb}{0.0,0.0,0.0}
\definecolor{cv5}{rgb}{0.0,0.0,0.0}
\definecolor{cfv5}{rgb}{1.0,1.0,1.0}
\definecolor{clv5}{rgb}{0.0,0.0,0.0}
\definecolor{cv6}{rgb}{0.0,0.0,0.0}
\definecolor{cfv6}{rgb}{1.0,1.0,1.0}
\definecolor{clv6}{rgb}{0.0,0.0,0.0}
\definecolor{cv7}{rgb}{0.0,0.0,0.0}
\definecolor{cfv7}{rgb}{1.0,1.0,1.0}
\definecolor{clv7}{rgb}{0.0,0.0,0.0}
\definecolor{cv8}{rgb}{0.0,0.0,0.0}
\definecolor{cfv8}{rgb}{1.0,1.0,1.0}
\definecolor{clv8}{rgb}{0.0,0.0,0.0}
\definecolor{cv9}{rgb}{0.0,0.0,0.0}
\definecolor{cfv9}{rgb}{1.0,1.0,1.0}
\definecolor{clv9}{rgb}{0.0,0.0,0.0}
\definecolor{cv0v2}{rgb}{0.0,0.0,0.0}
\definecolor{cv0v3}{rgb}{0.0,0.0,0.0}
\definecolor{cv0v7}{rgb}{0.0,0.0,0.0}
\definecolor{cv0v8}{rgb}{0.0,0.0,0.0}
\definecolor{cv1v3}{rgb}{0.0,0.0,0.0}
\definecolor{cv1v4}{rgb}{0.0,0.0,0.0}
\definecolor{cv1v8}{rgb}{0.0,0.0,0.0}
\definecolor{cv1v9}{rgb}{0.0,0.0,0.0}
\definecolor{cv2v4}{rgb}{0.0,0.0,0.0}
\definecolor{cv2v5}{rgb}{0.0,0.0,0.0}
\definecolor{cv2v9}{rgb}{0.0,0.0,0.0}
\definecolor{cv3v5}{rgb}{0.0,0.0,0.0}
\definecolor{cv3v6}{rgb}{0.0,0.0,0.0}
\definecolor{cv4v6}{rgb}{0.0,0.0,0.0}
\definecolor{cv4v7}{rgb}{0.0,0.0,0.0}
\definecolor{cv5v7}{rgb}{0.0,0.0,0.0}
\definecolor{cv5v8}{rgb}{0.0,0.0,0.0}
\definecolor{cv6v8}{rgb}{0.0,0.0,0.0}
\definecolor{cv6v9}{rgb}{0.0,0.0,0.0}
\definecolor{cv7v9}{rgb}{0.0,0.0,0.0}
\Vertex[style={minimum size=1.0cm,draw=cv0,fill=cfv0,text=clv0,shape=circle},LabelOut=false,L=\hbox{$0$},x=5.0cm,y=2.5cm]{v0}
\Vertex[style={minimum size=1.0cm,draw=cv1,fill=cfv1,text=clv1,shape=circle},LabelOut=false,L=\hbox{$1$},x=4.5225cm,y=4.0451cm]{v1}
\Vertex[style={minimum size=1.0cm,draw=cv2,fill=cfv2,text=clv2,shape=circle},LabelOut=false,L=\hbox{$2$},x=3.2725cm,y=5.0cm]{v2}
\Vertex[style={minimum size=1.0cm,draw=cv3,fill=cfv3,text=clv3,shape=circle},LabelOut=false,L=\hbox{$3$},x=1.7275cm,y=5.0cm]{v3}
\Vertex[style={minimum size=1.0cm,draw=cv4,fill=cfv4,text=clv4,shape=circle},LabelOut=false,L=\hbox{$4$},x=0.4775cm,y=4.0451cm]{v4}
\Vertex[style={minimum size=1.0cm,draw=cv5,fill=cfv5,text=clv5,shape=circle},LabelOut=false,L=\hbox{$5$},x=0.0cm,y=2.5cm]{v5}
\Vertex[style={minimum size=1.0cm,draw=cv6,fill=cfv6,text=clv6,shape=circle},LabelOut=false,L=\hbox{$6$},x=0.4775cm,y=0.9549cm]{v6}
\Vertex[style={minimum size=1.0cm,draw=cv7,fill=cfv7,text=clv7,shape=circle},LabelOut=false,L=\hbox{$7$},x=1.7275cm,y=0.0cm]{v7}
\Vertex[style={minimum size=1.0cm,draw=cv8,fill=cfv8,text=clv8,shape=circle},LabelOut=false,L=\hbox{$8$},x=3.2725cm,y=0.0cm]{v8}
\Vertex[style={minimum size=1.0cm,draw=cv9,fill=cfv9,text=clv9,shape=circle},LabelOut=false,L=\hbox{$9$},x=4.5225cm,y=0.9549cm]{v9}
\Edge[lw=0.02cm,style={color=cv0v2,},](v0)(v2)
\Edge[lw=0.02cm,style={color=cv0v3,},](v0)(v3)
\Edge[lw=0.02cm,style={color=cv0v7,},](v0)(v7)
\Edge[lw=0.02cm,style={color=cv0v8,},](v0)(v8)
\Edge[lw=0.02cm,style={color=cv1v3,},](v1)(v3)
\Edge[lw=0.02cm,style={color=cv1v4,},](v1)(v4)
\Edge[lw=0.02cm,style={color=cv1v8,},](v1)(v8)
\Edge[lw=0.02cm,style={color=cv1v9,},](v1)(v9)
\Edge[lw=0.02cm,style={color=cv2v4,},](v2)(v4)
\Edge[lw=0.02cm,style={color=cv2v5,},](v2)(v5)
\Edge[lw=0.02cm,style={color=cv2v9,},](v2)(v9)
\Edge[lw=0.02cm,style={color=cv3v5,},](v3)(v5)
\Edge[lw=0.02cm,style={color=cv3v6,},](v3)(v6)
\Edge[lw=0.02cm,style={color=cv4v6,},](v4)(v6)
\Edge[lw=0.02cm,style={color=cv4v7,},](v4)(v7)
\Edge[lw=0.02cm,style={color=cv5v7,},](v5)(v7)
\Edge[lw=0.02cm,style={color=cv5v8,},](v5)(v8)
\Edge[lw=0.02cm,style={color=cv6v8,},](v6)(v8)
\Edge[lw=0.02cm,style={color=cv6v9,},](v6)(v9)
\Edge[lw=0.02cm,style={color=cv7v9,},](v7)(v9)
\end{tikzpicture}
	\captionof{figure}{$X(\ints_{10},\{2,3,7,8\})$}
	\label{X10_23}
\end{minipage}%
\begin{minipage}{0.5\textwidth}
\centering
\begin{tikzpicture}
\definecolor{cv0}{rgb}{0.0,0.0,0.0}
\definecolor{cfv0}{rgb}{1.0,1.0,1.0}
\definecolor{clv0}{rgb}{0.0,0.0,0.0}
\definecolor{cv1}{rgb}{0.0,0.0,0.0}
\definecolor{cfv1}{rgb}{1.0,1.0,1.0}
\definecolor{clv1}{rgb}{0.0,0.0,0.0}
\definecolor{cv2}{rgb}{0.0,0.0,0.0}
\definecolor{cfv2}{rgb}{1.0,1.0,1.0}
\definecolor{clv2}{rgb}{0.0,0.0,0.0}
\definecolor{cv3}{rgb}{0.0,0.0,0.0}
\definecolor{cfv3}{rgb}{1.0,1.0,1.0}
\definecolor{clv3}{rgb}{0.0,0.0,0.0}
\definecolor{cv4}{rgb}{0.0,0.0,0.0}
\definecolor{cfv4}{rgb}{1.0,1.0,1.0}
\definecolor{clv4}{rgb}{0.0,0.0,0.0}
\definecolor{cv5}{rgb}{0.0,0.0,0.0}
\definecolor{cfv5}{rgb}{1.0,1.0,1.0}
\definecolor{clv5}{rgb}{0.0,0.0,0.0}
\definecolor{cv6}{rgb}{0.0,0.0,0.0}
\definecolor{cfv6}{rgb}{1.0,1.0,1.0}
\definecolor{clv6}{rgb}{0.0,0.0,0.0}
\definecolor{cv7}{rgb}{0.0,0.0,0.0}
\definecolor{cfv7}{rgb}{1.0,1.0,1.0}
\definecolor{clv7}{rgb}{0.0,0.0,0.0}
\definecolor{cv8}{rgb}{0.0,0.0,0.0}
\definecolor{cfv8}{rgb}{1.0,1.0,1.0}
\definecolor{clv8}{rgb}{0.0,0.0,0.0}
\definecolor{cv9}{rgb}{0.0,0.0,0.0}
\definecolor{cfv9}{rgb}{1.0,1.0,1.0}
\definecolor{clv9}{rgb}{0.0,0.0,0.0}
\definecolor{cv10}{rgb}{0.0,0.0,0.0}
\definecolor{cfv10}{rgb}{1.0,1.0,1.0}
\definecolor{clv10}{rgb}{0.0,0.0,0.0}
\definecolor{cv11}{rgb}{0.0,0.0,0.0}
\definecolor{cfv11}{rgb}{1.0,1.0,1.0}
\definecolor{clv11}{rgb}{0.0,0.0,0.0}
\definecolor{cv12}{rgb}{0.0,0.0,0.0}
\definecolor{cfv12}{rgb}{1.0,1.0,1.0}
\definecolor{clv12}{rgb}{0.0,0.0,0.0}
\definecolor{cv13}{rgb}{0.0,0.0,0.0}
\definecolor{cfv13}{rgb}{1.0,1.0,1.0}
\definecolor{clv13}{rgb}{0.0,0.0,0.0}
\definecolor{cv0v1}{rgb}{0.0,0.0,0.0}
\definecolor{cv0v6}{rgb}{0.0,0.0,0.0}
\definecolor{cv0v8}{rgb}{0.0,0.0,0.0}
\definecolor{cv0v13}{rgb}{0.0,0.0,0.0}
\definecolor{cv1v2}{rgb}{0.0,0.0,0.0}
\definecolor{cv1v7}{rgb}{0.0,0.0,0.0}
\definecolor{cv1v9}{rgb}{0.0,0.0,0.0}
\definecolor{cv2v3}{rgb}{0.0,0.0,0.0}
\definecolor{cv2v8}{rgb}{0.0,0.0,0.0}
\definecolor{cv2v10}{rgb}{0.0,0.0,0.0}
\definecolor{cv3v4}{rgb}{0.0,0.0,0.0}
\definecolor{cv3v9}{rgb}{0.0,0.0,0.0}
\definecolor{cv3v11}{rgb}{0.0,0.0,0.0}
\definecolor{cv4v5}{rgb}{0.0,0.0,0.0}
\definecolor{cv4v10}{rgb}{0.0,0.0,0.0}
\definecolor{cv4v12}{rgb}{0.0,0.0,0.0}
\definecolor{cv5v6}{rgb}{0.0,0.0,0.0}
\definecolor{cv5v11}{rgb}{0.0,0.0,0.0}
\definecolor{cv5v13}{rgb}{0.0,0.0,0.0}
\definecolor{cv6v7}{rgb}{0.0,0.0,0.0}
\definecolor{cv6v12}{rgb}{0.0,0.0,0.0}
\definecolor{cv7v8}{rgb}{0.0,0.0,0.0}
\definecolor{cv7v13}{rgb}{0.0,0.0,0.0}
\definecolor{cv8v9}{rgb}{0.0,0.0,0.0}
\definecolor{cv9v10}{rgb}{0.0,0.0,0.0}
\definecolor{cv10v11}{rgb}{0.0,0.0,0.0}
\definecolor{cv11v12}{rgb}{0.0,0.0,0.0}
\definecolor{cv12v13}{rgb}{0.0,0.0,0.0}
\Vertex[style={minimum size=1.0cm,draw=cv0,fill=cfv0,text=clv0,shape=circle},LabelOut=false,L=\hbox{$0$},x=5.0cm,y=2.5cm]{v0}
\Vertex[style={minimum size=1.0cm,draw=cv1,fill=cfv1,text=clv1,shape=circle},LabelOut=false,L=\hbox{$1$},x=4.7524cm,y=3.6126cm]{v1}
\Vertex[style={minimum size=1.0cm,draw=cv2,fill=cfv2,text=clv2,shape=circle},LabelOut=false,L=\hbox{$2$},x=4.0587cm,y=4.5048cm]{v2}
\Vertex[style={minimum size=1.0cm,draw=cv3,fill=cfv3,text=clv3,shape=circle},LabelOut=false,L=\hbox{$3$},x=3.0563cm,y=5.0cm]{v3}
\Vertex[style={minimum size=1.0cm,draw=cv4,fill=cfv4,text=clv4,shape=circle},LabelOut=false,L=\hbox{$4$},x=1.9437cm,y=5.0cm]{v4}
\Vertex[style={minimum size=1.0cm,draw=cv5,fill=cfv5,text=clv5,shape=circle},LabelOut=false,L=\hbox{$5$},x=0.9413cm,y=4.5048cm]{v5}
\Vertex[style={minimum size=1.0cm,draw=cv6,fill=cfv6,text=clv6,shape=circle},LabelOut=false,L=\hbox{$6$},x=0.2476cm,y=3.6126cm]{v6}
\Vertex[style={minimum size=1.0cm,draw=cv7,fill=cfv7,text=clv7,shape=circle},LabelOut=false,L=\hbox{$7$},x=0.0cm,y=2.5cm]{v7}
\Vertex[style={minimum size=1.0cm,draw=cv8,fill=cfv8,text=clv8,shape=circle},LabelOut=false,L=\hbox{$8$},x=0.2476cm,y=1.3874cm]{v8}
\Vertex[style={minimum size=1.0cm,draw=cv9,fill=cfv9,text=clv9,shape=circle},LabelOut=false,L=\hbox{$9$},x=0.9413cm,y=0.4952cm]{v9}
\Vertex[style={minimum size=1.0cm,draw=cv10,fill=cfv10,text=clv10,shape=circle},LabelOut=false,L=\hbox{$10$},x=1.9437cm,y=0.0cm]{v10}
\Vertex[style={minimum size=1.0cm,draw=cv11,fill=cfv11,text=clv11,shape=circle},LabelOut=false,L=\hbox{$11$},x=3.0563cm,y=0.0cm]{v11}
\Vertex[style={minimum size=1.0cm,draw=cv12,fill=cfv12,text=clv12,shape=circle},LabelOut=false,L=\hbox{$12$},x=4.0587cm,y=0.4952cm]{v12}
\Vertex[style={minimum size=1.0cm,draw=cv13,fill=cfv13,text=clv13,shape=circle},LabelOut=false,L=\hbox{$13$},x=4.7524cm,y=1.3874cm]{v13}
\Edge[lw=0.02cm,style={color=cv0v1,},](v0)(v1)
\Edge[lw=0.02cm,style={color=cv0v6,},](v0)(v6)
\Edge[lw=0.02cm,style={color=cv0v8,},](v0)(v8)
\Edge[lw=0.02cm,style={color=cv0v13,},](v0)(v13)
\Edge[lw=0.02cm,style={color=cv1v2,},](v1)(v2)
\Edge[lw=0.02cm,style={color=cv1v7,},](v1)(v7)
\Edge[lw=0.02cm,style={color=cv1v9,},](v1)(v9)
\Edge[lw=0.02cm,style={color=cv2v3,},](v2)(v3)
\Edge[lw=0.02cm,style={color=cv2v8,},](v2)(v8)
\Edge[lw=0.02cm,style={color=cv2v10,},](v2)(v10)
\Edge[lw=0.02cm,style={color=cv3v4,},](v3)(v4)
\Edge[lw=0.02cm,style={color=cv3v9,},](v3)(v9)
\Edge[lw=0.02cm,style={color=cv3v11,},](v3)(v11)
\Edge[lw=0.02cm,style={color=cv4v5,},](v4)(v5)
\Edge[lw=0.02cm,style={color=cv4v10,},](v4)(v10)
\Edge[lw=0.02cm,style={color=cv4v12,},](v4)(v12)
\Edge[lw=0.02cm,style={color=cv5v6,},](v5)(v6)
\Edge[lw=0.02cm,style={color=cv5v11,},](v5)(v11)
\Edge[lw=0.02cm,style={color=cv5v13,},](v5)(v13)
\Edge[lw=0.02cm,style={color=cv6v7,},](v6)(v7)
\Edge[lw=0.02cm,style={color=cv6v12,},](v6)(v12)
\Edge[lw=0.02cm,style={color=cv7v8,},](v7)(v8)
\Edge[lw=0.02cm,style={color=cv7v13,},](v7)(v13)
\Edge[lw=0.02cm,style={color=cv8v9,},](v8)(v9)
\Edge[lw=0.02cm,style={color=cv9v10,},](v9)(v10)
\Edge[lw=0.02cm,style={color=cv10v11,},](v10)(v11)
\Edge[lw=0.02cm,style={color=cv11v12,},](v11)(v12)
\Edge[lw=0.02cm,style={color=cv12v13,},](v12)(v13)
\end{tikzpicture}
	\captionof{figure}{$X(\ints_{14}, \{1,6,8,13\})$}
	\label{X14_16}
\end{minipage}
\begin{minipage}{0.5\textwidth}
\centering
\begin{tikzpicture}
\definecolor{cv0}{rgb}{0.0,0.0,0.0}
\definecolor{cfv0}{rgb}{1.0,1.0,1.0}
\definecolor{clv0}{rgb}{0.0,0.0,0.0}
\definecolor{cv1}{rgb}{0.0,0.0,0.0}
\definecolor{cfv1}{rgb}{1.0,1.0,1.0}
\definecolor{clv1}{rgb}{0.0,0.0,0.0}
\definecolor{cv2}{rgb}{0.0,0.0,0.0}
\definecolor{cfv2}{rgb}{1.0,1.0,1.0}
\definecolor{clv2}{rgb}{0.0,0.0,0.0}
\definecolor{cv3}{rgb}{0.0,0.0,0.0}
\definecolor{cfv3}{rgb}{1.0,1.0,1.0}
\definecolor{clv3}{rgb}{0.0,0.0,0.0}
\definecolor{cv4}{rgb}{0.0,0.0,0.0}
\definecolor{cfv4}{rgb}{1.0,1.0,1.0}
\definecolor{clv4}{rgb}{0.0,0.0,0.0}
\definecolor{cv5}{rgb}{0.0,0.0,0.0}
\definecolor{cfv5}{rgb}{1.0,1.0,1.0}
\definecolor{clv5}{rgb}{0.0,0.0,0.0}
\definecolor{cv6}{rgb}{0.0,0.0,0.0}
\definecolor{cfv6}{rgb}{1.0,1.0,1.0}
\definecolor{clv6}{rgb}{0.0,0.0,0.0}
\definecolor{cv7}{rgb}{0.0,0.0,0.0}
\definecolor{cfv7}{rgb}{1.0,1.0,1.0}
\definecolor{clv7}{rgb}{0.0,0.0,0.0}
\definecolor{cv8}{rgb}{0.0,0.0,0.0}
\definecolor{cfv8}{rgb}{1.0,1.0,1.0}
\definecolor{clv8}{rgb}{0.0,0.0,0.0}
\definecolor{cv9}{rgb}{0.0,0.0,0.0}
\definecolor{cfv9}{rgb}{1.0,1.0,1.0}
\definecolor{clv9}{rgb}{0.0,0.0,0.0}
\definecolor{cv10}{rgb}{0.0,0.0,0.0}
\definecolor{cfv10}{rgb}{1.0,1.0,1.0}
\definecolor{clv10}{rgb}{0.0,0.0,0.0}
\definecolor{cv11}{rgb}{0.0,0.0,0.0}
\definecolor{cfv11}{rgb}{1.0,1.0,1.0}
\definecolor{clv11}{rgb}{0.0,0.0,0.0}
\definecolor{cv12}{rgb}{0.0,0.0,0.0}
\definecolor{cfv12}{rgb}{1.0,1.0,1.0}
\definecolor{clv12}{rgb}{0.0,0.0,0.0}
\definecolor{cv13}{rgb}{0.0,0.0,0.0}
\definecolor{cfv13}{rgb}{1.0,1.0,1.0}
\definecolor{clv13}{rgb}{0.0,0.0,0.0}
\definecolor{cv0v2}{rgb}{0.0,0.0,0.0}
\definecolor{cv0v5}{rgb}{0.0,0.0,0.0}
\definecolor{cv0v9}{rgb}{0.0,0.0,0.0}
\definecolor{cv0v12}{rgb}{0.0,0.0,0.0}
\definecolor{cv1v3}{rgb}{0.0,0.0,0.0}
\definecolor{cv1v6}{rgb}{0.0,0.0,0.0}
\definecolor{cv1v10}{rgb}{0.0,0.0,0.0}
\definecolor{cv1v13}{rgb}{0.0,0.0,0.0}
\definecolor{cv2v4}{rgb}{0.0,0.0,0.0}
\definecolor{cv2v7}{rgb}{0.0,0.0,0.0}
\definecolor{cv2v11}{rgb}{0.0,0.0,0.0}
\definecolor{cv3v5}{rgb}{0.0,0.0,0.0}
\definecolor{cv3v8}{rgb}{0.0,0.0,0.0}
\definecolor{cv3v12}{rgb}{0.0,0.0,0.0}
\definecolor{cv4v6}{rgb}{0.0,0.0,0.0}
\definecolor{cv4v9}{rgb}{0.0,0.0,0.0}
\definecolor{cv4v13}{rgb}{0.0,0.0,0.0}
\definecolor{cv5v7}{rgb}{0.0,0.0,0.0}
\definecolor{cv5v10}{rgb}{0.0,0.0,0.0}
\definecolor{cv6v8}{rgb}{0.0,0.0,0.0}
\definecolor{cv6v11}{rgb}{0.0,0.0,0.0}
\definecolor{cv7v9}{rgb}{0.0,0.0,0.0}
\definecolor{cv7v12}{rgb}{0.0,0.0,0.0}
\definecolor{cv8v10}{rgb}{0.0,0.0,0.0}
\definecolor{cv8v13}{rgb}{0.0,0.0,0.0}
\definecolor{cv9v11}{rgb}{0.0,0.0,0.0}
\definecolor{cv10v12}{rgb}{0.0,0.0,0.0}
\definecolor{cv11v13}{rgb}{0.0,0.0,0.0}
\Vertex[style={minimum size=1.0cm,draw=cv0,fill=cfv0,text=clv0,shape=circle},LabelOut=false,L=\hbox{$0$},x=5.0cm,y=2.5cm]{v0}
\Vertex[style={minimum size=1.0cm,draw=cv1,fill=cfv1,text=clv1,shape=circle},LabelOut=false,L=\hbox{$1$},x=4.7524cm,y=3.6126cm]{v1}
\Vertex[style={minimum size=1.0cm,draw=cv2,fill=cfv2,text=clv2,shape=circle},LabelOut=false,L=\hbox{$2$},x=4.0587cm,y=4.5048cm]{v2}
\Vertex[style={minimum size=1.0cm,draw=cv3,fill=cfv3,text=clv3,shape=circle},LabelOut=false,L=\hbox{$3$},x=3.0563cm,y=5.0cm]{v3}
\Vertex[style={minimum size=1.0cm,draw=cv4,fill=cfv4,text=clv4,shape=circle},LabelOut=false,L=\hbox{$4$},x=1.9437cm,y=5.0cm]{v4}
\Vertex[style={minimum size=1.0cm,draw=cv5,fill=cfv5,text=clv5,shape=circle},LabelOut=false,L=\hbox{$5$},x=0.9413cm,y=4.5048cm]{v5}
\Vertex[style={minimum size=1.0cm,draw=cv6,fill=cfv6,text=clv6,shape=circle},LabelOut=false,L=\hbox{$6$},x=0.2476cm,y=3.6126cm]{v6}
\Vertex[style={minimum size=1.0cm,draw=cv7,fill=cfv7,text=clv7,shape=circle},LabelOut=false,L=\hbox{$7$},x=0.0cm,y=2.5cm]{v7}
\Vertex[style={minimum size=1.0cm,draw=cv8,fill=cfv8,text=clv8,shape=circle},LabelOut=false,L=\hbox{$8$},x=0.2476cm,y=1.3874cm]{v8}
\Vertex[style={minimum size=1.0cm,draw=cv9,fill=cfv9,text=clv9,shape=circle},LabelOut=false,L=\hbox{$9$},x=0.9413cm,y=0.4952cm]{v9}
\Vertex[style={minimum size=1.0cm,draw=cv10,fill=cfv10,text=clv10,shape=circle},LabelOut=false,L=\hbox{$10$},x=1.9437cm,y=0.0cm]{v10}
\Vertex[style={minimum size=1.0cm,draw=cv11,fill=cfv11,text=clv11,shape=circle},LabelOut=false,L=\hbox{$11$},x=3.0563cm,y=0.0cm]{v11}
\Vertex[style={minimum size=1.0cm,draw=cv12,fill=cfv12,text=clv12,shape=circle},LabelOut=false,L=\hbox{$12$},x=4.0587cm,y=0.4952cm]{v12}
\Vertex[style={minimum size=1.0cm,draw=cv13,fill=cfv13,text=clv13,shape=circle},LabelOut=false,L=\hbox{$13$},x=4.7524cm,y=1.3874cm]{v13}
\Edge[lw=0.02cm,style={color=cv0v2,},](v0)(v2)
\Edge[lw=0.02cm,style={color=cv0v5,},](v0)(v5)
\Edge[lw=0.02cm,style={color=cv0v9,},](v0)(v9)
\Edge[lw=0.02cm,style={color=cv0v12,},](v0)(v12)
\Edge[lw=0.02cm,style={color=cv1v3,},](v1)(v3)
\Edge[lw=0.02cm,style={color=cv1v6,},](v1)(v6)
\Edge[lw=0.02cm,style={color=cv1v10,},](v1)(v10)
\Edge[lw=0.02cm,style={color=cv1v13,},](v1)(v13)
\Edge[lw=0.02cm,style={color=cv2v4,},](v2)(v4)
\Edge[lw=0.02cm,style={color=cv2v7,},](v2)(v7)
\Edge[lw=0.02cm,style={color=cv2v11,},](v2)(v11)
\Edge[lw=0.02cm,style={color=cv3v5,},](v3)(v5)
\Edge[lw=0.02cm,style={color=cv3v8,},](v3)(v8)
\Edge[lw=0.02cm,style={color=cv3v12,},](v3)(v12)
\Edge[lw=0.02cm,style={color=cv4v6,},](v4)(v6)
\Edge[lw=0.02cm,style={color=cv4v9,},](v4)(v9)
\Edge[lw=0.02cm,style={color=cv4v13,},](v4)(v13)
\Edge[lw=0.02cm,style={color=cv5v7,},](v5)(v7)
\Edge[lw=0.02cm,style={color=cv5v10,},](v5)(v10)
\Edge[lw=0.02cm,style={color=cv6v8,},](v6)(v8)
\Edge[lw=0.02cm,style={color=cv6v11,},](v6)(v11)
\Edge[lw=0.02cm,style={color=cv7v9,},](v7)(v9)
\Edge[lw=0.02cm,style={color=cv7v12,},](v7)(v12)
\Edge[lw=0.02cm,style={color=cv8v10,},](v8)(v10)
\Edge[lw=0.02cm,style={color=cv8v13,},](v8)(v13)
\Edge[lw=0.02cm,style={color=cv9v11,},](v9)(v11)
\Edge[lw=0.02cm,style={color=cv10v12,},](v10)(v12)
\Edge[lw=0.02cm,style={color=cv11v13,},](v11)(v13)
\end{tikzpicture}
	\captionof{figure}{$X(\ints_{14}, \{2,5,9,12\})$}
	\label{X14_25}
\end{minipage}%
\begin{minipage}{0.5\textwidth}
\centering
\begin{tikzpicture}
\definecolor{cv0}{rgb}{0.0,0.0,0.0}
\definecolor{cfv0}{rgb}{1.0,1.0,1.0}
\definecolor{clv0}{rgb}{0.0,0.0,0.0}
\definecolor{cv1}{rgb}{0.0,0.0,0.0}
\definecolor{cfv1}{rgb}{1.0,1.0,1.0}
\definecolor{clv1}{rgb}{0.0,0.0,0.0}
\definecolor{cv2}{rgb}{0.0,0.0,0.0}
\definecolor{cfv2}{rgb}{1.0,1.0,1.0}
\definecolor{clv2}{rgb}{0.0,0.0,0.0}
\definecolor{cv3}{rgb}{0.0,0.0,0.0}
\definecolor{cfv3}{rgb}{1.0,1.0,1.0}
\definecolor{clv3}{rgb}{0.0,0.0,0.0}
\definecolor{cv4}{rgb}{0.0,0.0,0.0}
\definecolor{cfv4}{rgb}{1.0,1.0,1.0}
\definecolor{clv4}{rgb}{0.0,0.0,0.0}
\definecolor{cv5}{rgb}{0.0,0.0,0.0}
\definecolor{cfv5}{rgb}{1.0,1.0,1.0}
\definecolor{clv5}{rgb}{0.0,0.0,0.0}
\definecolor{cv6}{rgb}{0.0,0.0,0.0}
\definecolor{cfv6}{rgb}{1.0,1.0,1.0}
\definecolor{clv6}{rgb}{0.0,0.0,0.0}
\definecolor{cv7}{rgb}{0.0,0.0,0.0}
\definecolor{cfv7}{rgb}{1.0,1.0,1.0}
\definecolor{clv7}{rgb}{0.0,0.0,0.0}
\definecolor{cv8}{rgb}{0.0,0.0,0.0}
\definecolor{cfv8}{rgb}{1.0,1.0,1.0}
\definecolor{clv8}{rgb}{0.0,0.0,0.0}
\definecolor{cv9}{rgb}{0.0,0.0,0.0}
\definecolor{cfv9}{rgb}{1.0,1.0,1.0}
\definecolor{clv9}{rgb}{0.0,0.0,0.0}
\definecolor{cv10}{rgb}{0.0,0.0,0.0}
\definecolor{cfv10}{rgb}{1.0,1.0,1.0}
\definecolor{clv10}{rgb}{0.0,0.0,0.0}
\definecolor{cv11}{rgb}{0.0,0.0,0.0}
\definecolor{cfv11}{rgb}{1.0,1.0,1.0}
\definecolor{clv11}{rgb}{0.0,0.0,0.0}
\definecolor{cv12}{rgb}{0.0,0.0,0.0}
\definecolor{cfv12}{rgb}{1.0,1.0,1.0}
\definecolor{clv12}{rgb}{0.0,0.0,0.0}
\definecolor{cv13}{rgb}{0.0,0.0,0.0}
\definecolor{cfv13}{rgb}{1.0,1.0,1.0}
\definecolor{clv13}{rgb}{0.0,0.0,0.0}
\definecolor{cv0v3}{rgb}{0.0,0.0,0.0}
\definecolor{cv0v4}{rgb}{0.0,0.0,0.0}
\definecolor{cv0v10}{rgb}{0.0,0.0,0.0}
\definecolor{cv0v11}{rgb}{0.0,0.0,0.0}
\definecolor{cv1v4}{rgb}{0.0,0.0,0.0}
\definecolor{cv1v5}{rgb}{0.0,0.0,0.0}
\definecolor{cv1v11}{rgb}{0.0,0.0,0.0}
\definecolor{cv1v12}{rgb}{0.0,0.0,0.0}
\definecolor{cv2v5}{rgb}{0.0,0.0,0.0}
\definecolor{cv2v6}{rgb}{0.0,0.0,0.0}
\definecolor{cv2v12}{rgb}{0.0,0.0,0.0}
\definecolor{cv2v13}{rgb}{0.0,0.0,0.0}
\definecolor{cv3v6}{rgb}{0.0,0.0,0.0}
\definecolor{cv3v7}{rgb}{0.0,0.0,0.0}
\definecolor{cv3v13}{rgb}{0.0,0.0,0.0}
\definecolor{cv4v7}{rgb}{0.0,0.0,0.0}
\definecolor{cv4v8}{rgb}{0.0,0.0,0.0}
\definecolor{cv5v8}{rgb}{0.0,0.0,0.0}
\definecolor{cv5v9}{rgb}{0.0,0.0,0.0}
\definecolor{cv6v9}{rgb}{0.0,0.0,0.0}
\definecolor{cv6v10}{rgb}{0.0,0.0,0.0}
\definecolor{cv7v10}{rgb}{0.0,0.0,0.0}
\definecolor{cv7v11}{rgb}{0.0,0.0,0.0}
\definecolor{cv8v11}{rgb}{0.0,0.0,0.0}
\definecolor{cv8v12}{rgb}{0.0,0.0,0.0}
\definecolor{cv9v12}{rgb}{0.0,0.0,0.0}
\definecolor{cv9v13}{rgb}{0.0,0.0,0.0}
\definecolor{cv10v13}{rgb}{0.0,0.0,0.0}
\Vertex[style={minimum size=1.0cm,draw=cv0,fill=cfv0,text=clv0,shape=circle},LabelOut=false,L=\hbox{$0$},x=5.0cm,y=2.5cm]{v0}
\Vertex[style={minimum size=1.0cm,draw=cv1,fill=cfv1,text=clv1,shape=circle},LabelOut=false,L=\hbox{$1$},x=4.7524cm,y=3.6126cm]{v1}
\Vertex[style={minimum size=1.0cm,draw=cv2,fill=cfv2,text=clv2,shape=circle},LabelOut=false,L=\hbox{$2$},x=4.0587cm,y=4.5048cm]{v2}
\Vertex[style={minimum size=1.0cm,draw=cv3,fill=cfv3,text=clv3,shape=circle},LabelOut=false,L=\hbox{$3$},x=3.0563cm,y=5.0cm]{v3}
\Vertex[style={minimum size=1.0cm,draw=cv4,fill=cfv4,text=clv4,shape=circle},LabelOut=false,L=\hbox{$4$},x=1.9437cm,y=5.0cm]{v4}
\Vertex[style={minimum size=1.0cm,draw=cv5,fill=cfv5,text=clv5,shape=circle},LabelOut=false,L=\hbox{$5$},x=0.9413cm,y=4.5048cm]{v5}
\Vertex[style={minimum size=1.0cm,draw=cv6,fill=cfv6,text=clv6,shape=circle},LabelOut=false,L=\hbox{$6$},x=0.2476cm,y=3.6126cm]{v6}
\Vertex[style={minimum size=1.0cm,draw=cv7,fill=cfv7,text=clv7,shape=circle},LabelOut=false,L=\hbox{$7$},x=0.0cm,y=2.5cm]{v7}
\Vertex[style={minimum size=1.0cm,draw=cv8,fill=cfv8,text=clv8,shape=circle},LabelOut=false,L=\hbox{$8$},x=0.2476cm,y=1.3874cm]{v8}
\Vertex[style={minimum size=1.0cm,draw=cv9,fill=cfv9,text=clv9,shape=circle},LabelOut=false,L=\hbox{$9$},x=0.9413cm,y=0.4952cm]{v9}
\Vertex[style={minimum size=1.0cm,draw=cv10,fill=cfv10,text=clv10,shape=circle},LabelOut=false,L=\hbox{$10$},x=1.9437cm,y=0.0cm]{v10}
\Vertex[style={minimum size=1.0cm,draw=cv11,fill=cfv11,text=clv11,shape=circle},LabelOut=false,L=\hbox{$11$},x=3.0563cm,y=0.0cm]{v11}
\Vertex[style={minimum size=1.0cm,draw=cv12,fill=cfv12,text=clv12,shape=circle},LabelOut=false,L=\hbox{$12$},x=4.0587cm,y=0.4952cm]{v12}
\Vertex[style={minimum size=1.0cm,draw=cv13,fill=cfv13,text=clv13,shape=circle},LabelOut=false,L=\hbox{$13$},x=4.7524cm,y=1.3874cm]{v13}
\Edge[lw=0.02cm,style={color=cv0v3,},](v0)(v3)
\Edge[lw=0.02cm,style={color=cv0v4,},](v0)(v4)
\Edge[lw=0.02cm,style={color=cv0v10,},](v0)(v10)
\Edge[lw=0.02cm,style={color=cv0v11,},](v0)(v11)
\Edge[lw=0.02cm,style={color=cv1v4,},](v1)(v4)
\Edge[lw=0.02cm,style={color=cv1v5,},](v1)(v5)
\Edge[lw=0.02cm,style={color=cv1v11,},](v1)(v11)
\Edge[lw=0.02cm,style={color=cv1v12,},](v1)(v12)
\Edge[lw=0.02cm,style={color=cv2v5,},](v2)(v5)
\Edge[lw=0.02cm,style={color=cv2v6,},](v2)(v6)
\Edge[lw=0.02cm,style={color=cv2v12,},](v2)(v12)
\Edge[lw=0.02cm,style={color=cv2v13,},](v2)(v13)
\Edge[lw=0.02cm,style={color=cv3v6,},](v3)(v6)
\Edge[lw=0.02cm,style={color=cv3v7,},](v3)(v7)
\Edge[lw=0.02cm,style={color=cv3v13,},](v3)(v13)
\Edge[lw=0.02cm,style={color=cv4v7,},](v4)(v7)
\Edge[lw=0.02cm,style={color=cv4v8,},](v4)(v8)
\Edge[lw=0.02cm,style={color=cv5v8,},](v5)(v8)
\Edge[lw=0.02cm,style={color=cv5v9,},](v5)(v9)
\Edge[lw=0.02cm,style={color=cv6v9,},](v6)(v9)
\Edge[lw=0.02cm,style={color=cv6v10,},](v6)(v10)
\Edge[lw=0.02cm,style={color=cv7v10,},](v7)(v10)
\Edge[lw=0.02cm,style={color=cv7v11,},](v7)(v11)
\Edge[lw=0.02cm,style={color=cv8v11,},](v8)(v11)
\Edge[lw=0.02cm,style={color=cv8v12,},](v8)(v12)
\Edge[lw=0.02cm,style={color=cv9v12,},](v9)(v12)
\Edge[lw=0.02cm,style={color=cv9v13,},](v9)(v13)
\Edge[lw=0.02cm,style={color=cv10v13,},](v10)(v13)
\end{tikzpicture}
	\captionof{figure}{$X(\ints_{14},\{3,4,10,11\})$}
	\label{X14_34}
\end{minipage}
\end{figure}

\section{Open Problems}
In this paper, we gave a characterization of perfect state transfer in arc-reversal Grover walks, and constructed an infinite family of circulant graphs that exhibit this phenomenon. However, many questions still remain open; we list a few, in the hope of finding more examples of state transfer.

Since perfect state transfer at step $k$ implies periodicity at step $2k$, the first question we could ask is the following.
\begin{enumerate}[(i)]
\item Which regular graphs have periodic vertices?
\end{enumerate}
Theorem \ref{Thm_evsp} gives a characterization for periodic vertices. Although this is a local condition on the eigenvalue support of a vertex, it is satisfied when the entire graph is periodic, that is, when all eigenvalues of the graph are $d$ times the real parts of some $k$-th roots of unity. Hence, one might want to study graphs for which
\[z^n\phi\left(\frac{d}{2}\left(z + \frac{1}{z}\right)\right)\]
is a factor of $x^k-1$. In \cite{Yoshie2018}, Yoshie investigated periodic arc-reversal Grover walks on distance regular graphs, and found all Hamming graphs and Johnson graphs that are periodic.

Looking back at our definition of perfect state transfer, we see that it is because the initial state lives in $\col(D_t^T)$ that perfect state transfer can be characterized using graph spectra. In principle, for any unit vector $x$,
\[e_u\otimes x\]
could serve as the initial state that concentrates on $u$. Thus, the second question is to understand what happens if we relax the assumption on the initial state.
\begin{enumerate}[(i)]
\setcounter{enumi}{1}
\item Is there an example of perfect state transfer, where the initial state does not lie in $\col(D_t^T)$? If so, can we characterize such perfect state transfer purely in terms of the spectral decomposition of $X$?
\end{enumerate}

Finally, perfect state transfer is a rare phenomenon, partly because the eigenvalue support needs to satisfy a very strong condition. However, if we are interested in \textsl{pretty good state transfer from $u$ to $v$}, that is, for any $\epsilon>0$, there is a time $k$ at which the probability of being on $v$ is $\epsilon$-close to $1$, given initial state $e_u\otimes x$, then the eigenvalue conditions in Theorem \ref{Thm_pst} can be relaxed. Pretty good state transfer has been found in continuous quantum walks, see for example \cite{Banchi2017,Kempton2017,Coutinho2017a,Eisenberg2018}. We would like to see discrete analogues.

\begin{enumerate}[(i)]
\setcounter{enumi}{2}
\item Can we characterize pretty good state transfer in discrete quantum walks? Is there any graph with pretty good state transfer but without perfect state transfer?
\end{enumerate}

\section*{Acknowledgement}
The author would like to thank Ada Chan, Gabriel Coutinho, Chris Godsil, Krystal Guo, and Christino Tamon for their helpful discussion and generous comments.

\bibliographystyle{amsplain}
\bibliography{dqw.bib}

\providecommand{\bysame}{\leavevmode\hbox to3em{\hrulefill}\thinspace}
\providecommand{\MR}{\relax\ifhmode\unskip\space\fi MR }
\providecommand{\MRhref}[2]{%
  \href{http://www.ams.org/mathscinet-getitem?mr=#1}{#2}
}
\providecommand{\href}[2]{#2}
\begin{thebibliography}{10}

\bibitem{Aharonov2000}
Dorit Aharonov, Andris Ambainis, Julia Kempe, and Umesh Vazirani,
  \emph{{Quantum walks on graphs}}, ACM Press (2000), 50--59.

\bibitem{Angeles-Canul2009}
Ricardo Angeles-Canul, Rachael Norton, Michael Opperman, Christopher Paribello,
  Matthew Russell, and Christino Tamon, \emph{{On quantum perfect state
  transfer in weighted join graphs}}, International Journal of Quantum
  Information \textbf{7} (2009), no.~8, 1429--1445.

\bibitem{Angeles-Canul2010}
\bysame, \emph{{Perfect state transfer, integral circulants and join of
  graphs}}, Quantum Information and Computation (2010), 325--342.

\bibitem{Bachman2012}
Rachel Bachman, Eric Fredette, Jessica Fuller, Michael Landry, Michael
  Opperman, Christino Tamon, and Andrew Tollefson, \emph{{Perfect state
  transfer of quantum walks on quotient graphs}}, Quantum Information and
  Computation (2012), 293--313.

\bibitem{Banchi2017}
Leonardo Banchi, Gabriel Coutinho, Chris Godsil, and Simone Severini,
  \emph{{Pretty good state transfer in qubit chains—The Heisenberg
  Hamiltonian}}, Journal of Mathematical Physics \textbf{58} (2017), no.~3.

\bibitem{Barr2014}
K.~Barr, T.~Proctor, D.~Allen, and V.~Kendon, \emph{{Periodicity and perfect
  state transfer in quantum walks on variants of cycles}}, Quantum Information
  {\&} Computation (2014), 417--438.

\bibitem{Bose2002}
Sougato Bose, \emph{{Quantum Communication through an unmodulated Spin Chain}},
  Physical Review Letters (2003).

\bibitem{Cheung2011}
Wang~Chi Cheung and Chris Godsil, \emph{{Perfect state transfer in cubelike
  graphs}}, Linear Algebra and Its Applications (2011), 2468--2474.

\bibitem{Childs2009a}
Andrew~M. Childs, \emph{{Universal computation by quantum walk}}, Physical
  Review Letters \textbf{102} (2009), no.~18.

\bibitem{Christandl2005}
Matthias Christandl, Nilanjana Datta, Tony~C. Dorlas, Artur Ekert, Alastair
  Kay, and Andrew~J. Landahl, \emph{{Perfect transfer of arbitrary states in
  quantum spin networks}}, Physical Review A - Atomic, Molecular, and Optical
  Physics \textbf{71} (2005), no.~3.

\bibitem{Christandl2004}
Matthias Christandl, Nilanjana Datta, Artur Ekert, and Andrew~J. Landahl,
  \emph{{Perfect state transfer in quantum spin networks}}, Physical Review
  Letters \textbf{92} (2004), no.~18.

\bibitem{Coutinho2014}
Gabriel Coutinho, \emph{{Quantum State Transfer in Graphs}}, Ph.D. thesis,
  University of Waterloo, 2014.

\bibitem{Coutinho2015a}
Gabriel Coutinho, Chris Godsil, Krystal Guo, and Fr{\'{e}}d{\'{e}}ric Vanhove,
  \emph{{Perfect state transfer on distance-regular graphs and association
  schemes}}, Linear Algebra and Its Applications \textbf{478} (2015), 108--130.

\bibitem{Coutinho2015}
Gabriel Coutinho and Chris~D Godsil, \emph{{Perfect state transfer in products
  and covers of graphs}}, Linear and Multilinear Algebra (2015), 1--12.

\bibitem{Coutinho2017a}
Gabriel Coutinho, Krystal Guo, and Christopher~M. {Van Bommel}, \emph{{Pretty
  good state transfer between internal nodes of paths}}, vol.~17, Rinton Press,
  2017.

\bibitem{Coutinho2015b}
Gabriel Coutinho and Henry Liu, \emph{{No Laplacian perfect state transfer in
  trees}}, SIAM Journal on Discrete Mathematics \textbf{29} (2015), no.~4,
  2179--2188.

\bibitem{Eisenberg2018}
Or~Eisenberg, Mark Kempton, and Gabor Lippner, \emph{{Pretty good quantum state
  transfer in asymmetric graphs via potential}}, arXiv:1804.01645 (2018).

\bibitem{Godsil2011a}
Chris Godsil, \emph{{Periodic graphs}}, The Electronic Journal of Combinatorics
  \textbf{18} (2011), no.~1.

\bibitem{Godsil2015a}
\bysame, \emph{{Graph Spectra and Quantum Walks}}, Unpublished, 2015.

\bibitem{Godsil2001}
Chris Godsil and Gordon Royle, \emph{{Algebraic Graph Theory}}, Springer New
  York, 2001.

\bibitem{Godsil2017Smith}
Chris Godsil and Jamie Smith, \emph{{Strongly cospectral vertices}},
  arXiv:1709.07975 (2017).

\bibitem{Kay2006}
Alastair Kay, \emph{{Perfect state transfer: Beyond nearest-neighbor
  couplings}}, Physical Review A - Atomic, Molecular, and Optical Physics
  \textbf{73} (2006), no.~3.

\bibitem{Kay2011}
\bysame, \emph{{Basics of perfect communication through quantum networks}},
  Physical Review A - Atomic, Molecular, and Optical Physics \textbf{84}
  (2011), no.~2.

\bibitem{Kempton2017}
Mark Kempton, Gabor Lippner, and Shing-Tung Yau, \emph{{Pretty good quantum
  state transfer in symmetric spin networks via magnetic field}}, Quantum
  Information Processing \textbf{16} (2017), no.~9.

\bibitem{Kendon2003}
Viv Kendon, \emph{{Quantum walks on general graphs}}, International Journal of
  Quantum Information (2003), 791--805.

\bibitem{Kendon2011}
Vivien~M. Kendon and Christino Tamon, \emph{{Perfect state transfer in quantum
  walks on graphs}}, 2011, pp.~422--433.

\bibitem{Kurzynski2011}
Pawe{\l} Kurzy{\'{n}}ski and Antoni W{\'{o}}jcik, \emph{{Discrete-time quantum
  walk approach to state transfer}}, Physical Review A - Atomic, Molecular, and
  Optical Physics \textbf{83} (2011), no.~6.

\bibitem{Lovett2010}
Neil~B Lovett, Sally Cooper, Matthew Everitt, Matthew Trevers, and Viv Kendon,
  \emph{{Universal quantum computation using the discrete-time quantum walk}},
  Physical Review A \textbf{81} (2010), no.~4.

\bibitem{Portugal2015}
R.~Portugal, R.~A.~M. Santos, T.~D. Fernandes, and D.~N. Gon{\c{c}}alves,
  \emph{{The staggered quantum walk model}}, Quantum Information Processing
  \textbf{15} (2016), no.~1, 85--101.

\bibitem{Stefanak2016}
Martin Stefanak and Stanislav Skoupy, \emph{{Perfect state transfer by means of
  discrete-time quantum walk search algorithms on highly symmetric graphs}},
  Physical Review A - Atomic, Molecular, and Optical Physics \textbf{94}
  (2016), no.~2.

\bibitem{Stefanak2017}
\bysame, \emph{{Perfect state transfer by means of discrete-time quantum walk
  on complete bipartite graphs}}, Quantum Information Processing \textbf{16}
  (2017), no.~3.

\bibitem{Szegedy2004}
Mario Szegedy, \emph{{Quantum speed-up of Markov chain based algorithms}}, 45th
  Annual IEEE Symposium on Foundations of Computer Science (2004), 32--41.

\bibitem{Underwood2010}
Michael Underwood and David Feder, \emph{{Universal quantum computation by
  discontinuous quantum walk}}, Physical Review A \textbf{82} (2010), no.~4,
  1--69.

\bibitem{Yalcnkaya2014}
I~Yal{\c{c}}ınkaya and Z~Gedik, \emph{{Qubit state transfer via
  one-dimensional discrete-time quantum walk}}, arXiv:1407.0689 (2014).

\bibitem{Yoshie2018}
Yusuke Yoshie, \emph{{Periodicity of Grover walks on distance-regular graphs}},
  arXiv:1805.07681 (2018).

\bibitem{Zhan2014}
Xiang Zhan, Hao Qin, Zhi-hao Bian, Jian Li, and Peng Xue, \emph{{Perfect state
  transfer and efficient quantum routing: a discrete-time quantum walk
  approach}}, Physical Review A - Atomic, Molecular, and Optical Physics
  (2014).

\end{thebibliography}
\end{document}